\documentclass[11pt]{amsart}
\usepackage{amsmath, amsfonts, amssymb, amsthm, amscd}

\usepackage{color}
\usepackage{hyperref}
\usepackage{graphicx}
\usepackage{tikz}

\usepackage[leqno]{amsmath}
\makeatletter
\newcommand{\leqnomode}{\tagsleft@true}
\newcommand{\reqnomode}{\tagsleft@false}
\makeatother

\usepackage{cases}
\usepackage{subfigure}
\usepackage{color}
\usepackage{latexsym, longtable}
\usepackage{enumerate}
\usepackage{comment}

\usepackage[T1]{fontenc}
\usepackage{tikz-cd}

\newtheorem{theorem}{Theorem}[section]
\newtheorem{proposition}[theorem]{Proposition}

\newtheorem{corollary}[theorem]{Corollary}

\newtheorem{definition}[theorem]{Definition}

\newtheorem{example}[theorem]{Example}

\newtheorem{remark}[theorem]{Remark}
\numberwithin{equation}{section}

\DeclareMathOperator{\perm}{perm}

\DeclareMathOperator{\ST}{ST}

\newcommand{\Z}{\mathbb{Z}}

\DeclareMathOperator{\arm}{\text{arm}}
\DeclareMathOperator{\leg}{\text{leg}}
\DeclareMathOperator{\Des}{Des}

\DeclareMathOperator{\maj}{maj}

\DeclareMathOperator{\inv}{inv}
\DeclareMathOperator{\inc}{inc}
\DeclareMathOperator{\coinv}{coinv}

\DeclareMathOperator{\wt}{wt}

\newlength\cellsize 
\savebox2{%
\begin{picture}(12,12)
\put(0,0){\line(1,0){12}}
\put(0,0){\line(0,1){12}}
\put(12,0){\line(0,1){12}}
\put(0,12){\line(1,0){12}}
\end{picture}}
\newcommand\cellify[1]{\def\thearg{#1}\def\nothing{}%
\ifx\thearg\nothing
\vrule width0pt height\cellsize depth0pt\else
\hbox to 0pt{\usebox2\hss}\fi%
\vbox to 12\unitlength{
\vss
\hbox to 12\unitlength{\hss$#1$\hss}
\vss}}
\newcommand\tableau[1]{\vtop{\let\\=\cr
\setlength\baselineskip{-16000pt}
\setlength\lineskiplimit{16000pt}
\setlength\lineskip{0pt}
\halign{&\cellify{##}\cr#1\crcr}}}
\savebox3{%
\begin{picture}(12,12)
\put(0,0){\line(1,0){12}}
\put(0,0){\line(0,1){12}}
\put(12,0){\line(0,1){12}}
\put(0,12){\line(1,0){12}}
\end{picture}}
\newcommand\expath[1]{%
\hbox to 0pt{\usebox3\hss}%
\vbox to 12\unitlength{
\vss
\hbox to 12\unitlength{\hss$#1$\hss}
\vss}}

\newcommand\cell[3]{
\def\i{#1} \def\j{#2} \def\entry{#3}
\draw (\j-1,-\i)--(\j,-\i)--(\j,-\i+1)--(\j-1,-\i+1)--(\j-1,-\i);
\node at (\j-.5,-\i+.5) {\entry};
}

\let\OLDthebibliography\thebibliography
\renewcommand\thebibliography[1]{
  \OLDthebibliography{#1}
  \setlength{\parskip}{0pt}
  \setlength{\itemsep}{3pt plus 0.5ex}
}

\title[Compact formulas for  Macdonald polynomials]{Compact formulas for Macdonald polynomials \\ and quasisymmetric Macdonald polynomials}

\date{\today}
\author[Corteel]{Sylvie Corteel}
\address{Department of Mathematics, UC Berkeley, USA}
\email{corteel@berkeley.edu}
\author[Haglund]{Jim Haglund}
\address{Department of Mathematics, University of Pennsylvania, USA}
\email{haglund@math.upenn.edu}
\thanks{Jim Haglund was partially supported by NSF grant DMS-1600670.}
\author[Mandelshtam]{Olya Mandelshtam}
\address{Department of Mathematics, Brown University USA}
\email{olya@math.brown.edu}
\thanks{Olya Mandelshtam was partially supported by NSF grant DMS-1704874.}
\author[Mason]{Sarah Mason}
\address{Department of Mathematics, Wake Forest University, USA}
\email{masonsk@wfu.edu}
\author[Williams]{Lauren Williams}
\address{Department of Mathematics, Harvard University, USA}
\email{williams@math.harvard.edu}
\thanks{Lauren Williams was partially supported by NSF grant DMS-1854512.}

\begin{document}
\begin{abstract}
We present several new and compact formulas for the modified and integral form of the Macdonald polynomials, building on the compact ``multiline queue'' formula for Macdonald polynomials due to Corteel, Mandelshtam and Williams. We also introduce a new quasisymmetric analogue of Macdonald polynomials.  These \emph{quasisymmetric Macdonald polynomials} refine the (symmetric) Macdonald polynomials and specialize at $q=t=0$ to the quasisymmetric Schur polynomials defined by Haglund, Luoto, Mason, and van Willigenburg. This is an extended abstract.
\end{abstract}

\maketitle

\section{Introduction}

The symmetric 
\emph{Macdonald polynomials} $P_{\lambda}(X; q, t)$ \cite{Macdonald}
are a family of polynomials in $X = \{x_1, x_2,\dots \}$ indexed by partitions, whose coefficients depend on two parameters $q$ and $t$.  Macdonald polynomials generalize 
multiple important families of polynomials, including Schur polynomials and Hall-Littlewood polynomials. They can be defined as the unique monic basis for the ring of symmetric functions that satisfies certain triangularity and orthogonality conditions. The related \emph{nonsymmetric Macdonald polynomials} $E_{\mu}(X;q,t)$ \cite{Mac88, Macdonald, Cher1} were introduced shortly after the introduction of Macdonald polynomials as a tool to study Macdonald polynomials. The $E_{\mu}(X;q, t)$ are indexed by weak compositions and form a basis for the full polynomial ring $\mathbb{Q}[X](q,t)$. 

There has been a great deal of work devoted to understanding Macdonald polynomials from a combinatorial point of view. Haglund-Haiman-Loehr \cite{HHL04} gave a combinatorial formula for the \emph{integral forms} $J_{\lambda}(X; q,t)$, which are scalar multiples of the classical monic forms $P_{\lambda}(X;q, t)$.  They also gave a formula for the nonsymmetric Macdonald polynomials $E_{\mu}(X;q,t)$ \cite{HHL08}, and for the \emph{transformed} or \emph{modified} Macdonald polynomials $\widetilde{H}_{\lambda}(X;q,t)$, which are obtained from $J_{\lambda}(X;q, t)$ via \emph{plethysm}. Macdonald conjectured and Haiman proved \cite{Hai99}, using the geometry of the Hilbert scheme, that the modified Macdonald polynomials $\widetilde{H}_{\lambda}(X;q,t)$ have a positive Schur expansion whose coefficients are $qt$-Kostka polynomials.  However, it is still an open problem to give a combinatorial proof of the Schur positivity or a manifestly positive formula for the $qt$-Kostka polynomials.

Recently, a beautiful connection has been found between Macdonald polynomials and a statistical mechanics model called the multispecies \emph{asymmetric simple exclusion process} (ASEP) on a circle. The ASEP is a one-dimensional exactly solvable particle model; Cantini-deGier-Wheeler \cite{CGW} showed that the partition function of the multispecies ASEP on a circle
is equal to a Macdonald polynomial $P_{\lambda}(x_1,\dots,x_n; q, t)$ evaluated at $q=1$ and $x_i=1$ for all $i$. Building on this result as well as work of Martin \cite{Martin}, the first, third, and fifth authors recently used \emph{multiline queues} to simultaneously compute the stationary probabilities of the multispecies exclusion process, and give compact formulas for the symmetric Macdonald polynomials $P_{\lambda}$ and the nonsymmetric Macdonald polynomials  $E_{\lambda}$ \cite{CMW18}, for any partition $\lambda$. These formulas are ``compact'' in that they have fewer terms than the formulas of Haglund-Haiman-Loehr.

In this paper, we use the above ideas to continue the search for compact formulas for Macdonald polynomials. Our first two main results are compact formulas for the modified Macdonald polynomials $\widetilde{H}_{\lambda}(X;q,t)$ and the integral forms $J_{\lambda}(X;q,t)$; these new formulas have far fewer terms than other known combinatorial formulas. Our third main result uses the connection with the ASEP on a ring towards a different application: the introduction of a new family of quasisymmetric functions we call \emph{quasisymmetric Macdonald polynomials}  $G_{\gamma}(X; q, t)$. We show that $G_{\gamma}(X; q, t)$ is indeed a quasisymmetric function, and give a combinatorial formula for the $G_{\gamma}(X; q, t)$ corresponding to ``pieces'' of the compact formula for the $P_{\lambda}(X; q, t)$ from \cite{CMW18}. Moreover, the quasisymmetric function $G_{\gamma}(X; q, t)$ at $q=t=0$ specializes to the \emph{quasisymmetric Schur functions} $\text{QS}_{\gamma}(X)$ introduced by the second and fourth authors, together with Luoto and van Willigenburg \cite{HLMV09}.  The quasisymmetric Schur functions form a basis for the ring of quasisymmetric functions, and until now it has been an open question to find a refinement of the Macdonald polynomials $P_{\lambda}$ into quasisymmetric pieces which generalize the quasisymmetric Schur functions. 

This paper is organized as follows. In Section~\ref{sec:def}, we provide the relevant background. Sections~\ref{sec:Hcompact} and~\ref{sec:Jcompact} describe our two compact formulas, and Section~\ref{sec:qsym} defines our new quasisymmetric Macdonald polynomials.  While many open problems naturally arise from this work, we will defer their discussion to the longer version of this paper.

\section{Definitions}\label{sec:def}

We begin by introducing relevant notation and definitions.  In our partition and composition diagrams, given in French notation, the columns are labeled from left to right, and the rows are labeled from bottom to top, so that the notation $(i,r)$ refers to the box (or \emph{cell}) in the $i^{th}$ column from the left and the $r^{th}$ row from the bottom. Given a partition/composition $\alpha$, its diagram is a sequence of columns bottom justified, where the $i^{th}$ column has $\alpha_i$ cells.  The \emph{leg} of a cell $(i,r)$, denoted $\leg((i,r))$, equals the number of cells in column $i$ above the cell $(i,r)$.  Analogously the \emph{arm} of a cell $(i,r)$, denoted $\arm((i,r))$, equals the number of cells in row $r$ to the right of the cell $(i,r)$.   
 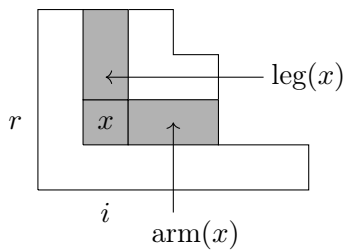
\begin{figure}[h]
 \begin{center}
\begin{tikzpicture}[scale=.6]
\draw (0,0)--(6,0)--(6,1)--(4,1)--(4,3)--(3,3)--(3,4)--(2,4)--(0,4)--(0,0);
\draw[fill=gray!60] (1,1)--(2,1)--(2,2)--(1,2)--(1,1);
\draw[fill=gray!60] (2,2)--(2,4)--(1,4)--(1,2)--(2,2);
\draw[fill=gray!60] (2,2)--(4,2)--(4,1)--(2,1)--(2,2);
\node at (1.5,1.5) {$x$};
\node at (1.5,-.5) {$i$};
\node at (-.5,1.5) {$r$};
\draw[->] (5,2.5)--(1.5,2.5);
\node at (6,2.5) {$\leg(x)$} ;
\draw[->] (3,-.5)--(3,1.5);
\node at (3.5,-1) {$\arm(x)$} ;
\end{tikzpicture}
\end{center}
 \caption{$\leg(x)$ of cell $x$ in row $r$ and column $i$.}
\label{fig:leg}
\end{figure}

A \emph{filling} $\sigma:\lambda \to \Z^+$ is an assignment of positive integers to the cells of $\lambda$ and is denoted by $(\sigma, \lambda)$.  For $s\in \lambda$, let $\sigma(s)$ denote the integer assigned to $s$, i.e., the integer occupying cell $s$.  The numbers appearing in such a filling are called the \emph{entries}. For each filling $\sigma$ of $\lambda$ we associate $x$, $q$, and $t$ \emph{weights}.  The \emph{$x$-weight} is defined in a similar fashion to semistandard Young tableaux, namely \[x^{\sigma} = \prod_{s\in \lambda} x_{\sigma(s)}.\]

We recall several definitions from~\cite{HHL04}. Assume that the diagram of a partition $\lambda$ has a \emph{basement}, i.e., a zero(th) row of size $\lambda_1$ all of whose cells are filled with the entry $\infty$. Let $u$, $v$, and $r$ be positive integers with $u<v$. Given a (diagram of a) partition $\lambda$ and a filling $(\sigma, \lambda)$, a \emph{triple} consists of the three cells (if they are present in the diagram) $(v,r)$, $(u,r)$, and $(u,r-1)$. Let $a = \sigma((v,r))$, $b = \sigma((u,r))$, and $c=\sigma((u,r-1))$. We say that the triple  is a \emph{counterclockwise inversion triple} if any of the following holds:
\begin{align*}
a < b \leq c\ \ {\rm or}\ \ 
c < a < b \ \ {\rm or}\ \ 
b \leq c < a.
\end{align*}
For example, in Figure \ref{fig:perm} the entries $(3,3)$, $(1,3)$ and $(1,4)$  form a counterclockwise inversion triple. We say that the triple  is a \emph{clockwise inversion triple} if any of the following holds:
\begin{align*}
a > b > c\ \ {\rm or}\ \
c > a > b \ \ {\rm or}\ \
b > c > a.
\end{align*}

Note that since $\sigma((j,0))=\infty$, if $r=1$ then $(v,1)$ and  $(u,1)$ form a (counterclockwise) inversion triple if and only if $\sigma((u,1))>\sigma((v,1))$. In this case we say that the triple is {\em degenerate}. For example, in Figure \ref{fig:perm} the entries $(3,1)$ and $(6,1)$  form a degenerate inversion triple. Given a filling $(\sigma, \lambda)$, let $\inv(\sigma, \lambda)$ be the total number of  (counterclockwise) inversion triples including degenerate triples.  Let $\coinv(\sigma, \lambda)$ be the total number of triples in the filling minus $\inv(\sigma, \lambda)$.  (In other words, $\coinv(\sigma, \lambda)$ is the number of triples which are not counterclockwise inversion triples.)

Define the \emph{set of descents} of a filling to be 
\[
\Des(\lambda, \sigma) = \{(u,r)\in \lambda\ :\ \sigma((u,r))>\sigma((u,r-1))\}.
\]   
We define the \emph{major index} $\maj(\sigma)$ to be the sum over the legs of the descents of 
	$(\lambda, \sigma)$:
\[
\maj(\sigma) = \sum_{x \in \Des(\sigma)} (\leg(x)+1).
\]
The tableau $\sigma$ of shape $\lambda$ in Figure \ref{fig:perm} has $\inv(\sigma,\lambda)=22$ and $\maj(\sigma,\lambda)=5$. Two cells are \emph{attacking} if their entries are equal and the cells are either in the same row, or they are in adjacent rows, with the rightmost cell in a row strictly below the other cell. A filling is \emph{nonattacking} if it does not contain any attacking pairs of cells.

Given a \emph{weak composition}, i.e., a vector $\alpha = (\alpha _1, \alpha _2, \ldots , \alpha _n)$ of nonnegative integers, we let $\inc(\alpha)$ and $\text{dec}(\alpha)$ be the vectors obtained from $\alpha$ by sorting the parts in weakly increasing order, and weakly decreasing order, respectively. Let $\beta (\alpha)$ be the permutation in $S_n$ of maximal length with the property that $\beta$ applied to the vector $\alpha$ yields $\inc(\alpha)$.  Let $\alpha ^{+}$ be the \emph{strong composition} obtained from $\alpha$ by removing  the zeros and let $\ell (\alpha)$ be the number of parts of $\alpha ^{+}$.  For example, if $\alpha = (0,2,0,2,1, 3)$ then  $\inc(\alpha) = (0,0,1,2,2,3)$, $\text{dec}(\alpha) = (3,2,2,1,0,0)$, $\beta (\alpha) = (3,1,5,4,2,6)$, $\alpha ^{+} = (2,2,1,3)$, and $\ell (\alpha) = 4$.

\section{Compact formula for modified Macdonald polynomials}\label{sec:Hcompact}

\begin{figure}

  \centering

\begin{tikzpicture}[scale=.5]

\node at (-1.5,.5) {\rm row 5};
\node at (-1.5,-.5) {\rm row 4};
\node at (-1.5,-1.5) {\rm row 3};
\node at (-1.5,-2.5) {\rm row 2};
\node at (-1.5,-3.5) {\rm row 1};
\node at (-1.5,-4.5) {\rm row 0};

\cell315 \cell325 \cell336 \cell346 \cell356

\cell215 \cell225 \cell232  

\cell419 \cell429 \cell439\cell441 \cell453 \cell462 \cell472 \cell483 \cell493

\cell111 \cell121\cell131

\cell016\cell026\cell036

\cell51{$\infty$} \cell52{$\infty$}  \cell53{$\infty$} \cell54{$\infty$}  \cell55{$\infty$}  \cell56{$\infty$}  \cell57{$\infty$}  \cell58{$\infty$}  \cell59{$\infty$}

\draw (0,-4)--(0,1)--(3,1)--(3,-2)--(5,-2)--(5,-3)--(9,-3);

\end{tikzpicture} 
\caption{A sorted tableau $\sigma$, with 
	$\perm_t(\sigma) = {3 \choose 2,1}_t {2 \choose 1,1}_t {4 \choose 2,2}_t$.}
\label{fig:perm}
\end{figure}

Our first main result is a ``compact'' formula for the modified Macdonald
polynomials  $\widetilde{H}_{\lambda}(X; q, t)$.  Before explaining
our result, we first
recall the combinatorial formula 
 of Haiman, Haglund and Loehr \cite{HHL04}.
\begin{theorem}[{\cite{HHL04}}] \label{thm:HHL}
        The modified Macdonald polynomial
        $\widetilde{H}_{\lambda}(X;q,t)$
        is given by
\[
        \widetilde{H}_{\lambda}(X;q,t) = \sum_{\sigma: \lambda \to \Z^+} x^\sigma q^{\inv(\sigma, \lambda)} t^{\maj(\sigma, \lambda)}.
\]
\end{theorem}

While Theorem \ref{thm:HHL} is simple and elegant, it has the disadvantage of having many terms, since it is a sum over \emph{all} fillings of a Young diagram by positive integers. By contrast, our compact formula (Theorem \ref{thm:MacPerm2}) is a sum over far fewer terms---it is a sum over \emph{sorted tableaux}.  To define these sorted tableaux, we first define an order on the columns of the fillings.
\begin{definition}\label{def:orderhhl}
Fix a filling $\sigma$ of $\lambda$ and consider two columns $A$ and $B$ of height $j$ in $\lambda$. Let $a_1,\dots, a_j$ and $b_1,\dots, b_j$ be the entries of columns $A$ and $B$, respectively, read from bottom to top. We say that $A \lhd B$, if either $a_1 < b_1$, or $a_i = b_i$ for $i=1, 2,  \dots, h-1$ (for some positive $h$), and the cells containing $b_h$, $a_h$ and $a_{h-1}$ do not form an inversion triple.
\end{definition}

\begin{definition}\label{def:ptableau2}
Given a filling $\sigma$ of the diagram of a partition $\lambda$, we say that $\sigma$ is a \emph{sorted tableau} if, for all positive integers $h$, when we read all columns of height $h$ from left to right, the columns appear in weakly increasing order with respect to $\lhd$.  We write $\ST(\lambda)$ for the set of all sorted tableaux of shape $\lambda$.
\end{definition}

Let $\sigma$ be a sorted tableau. First, suppose that the shape of $\lambda$ is an $m \times n$ rectangle.  The $n$ columns may not all have distinct fillings: suppose that among those $n$ columns, there are $j$ distinct column fillings, with $u_1$ identical columns of the first filling, $u_2$ identical columns of the second filling, up to $u_j$ identical columns of the $j$th filling.  Define 
\[
\perm_t(\sigma) = {n \choose{u_1,\ldots,u_j}}_t.
\]

Suppose $\sigma$ is a sorted tableau which is a concatenation of rectangular sorted tableaux $\sigma_1,\dots, \sigma_{\ell}$, all of different heights.  Define $\perm_t(\sigma) = \prod_{i=1}^{\ell} \perm_t(\sigma_i).$ See Figure \ref{fig:perm} for an example.

Our main result in this section is a  compact formula for $\widetilde{H}_{\lambda'}$.
\begin{theorem}{\label{thm:MacPerm2}}
The modified Macdonald polynomial
$\widetilde{H}_{\lambda'}(X;q,t)$
	equals 
	\[
	\widetilde{H}_{\lambda'}(X;q,t) = \sum_{\sigma\in \ST(\lambda)} x^\sigma t^{\inv(\sigma, \lambda)} q^{\maj(\sigma, \lambda)} 
	\perm_t(\sigma, \lambda),
\]
where the sum is over all sorted tableaux of shape $\lambda$.
\end{theorem}

\begin{example} We use Theorem \ref{thm:MacPerm2} to compute
 $\widetilde{H}_{2,1,1}(x_1, x_2, x_3;q,t)$
 in Figure \ref{fig:alternativeDef}.
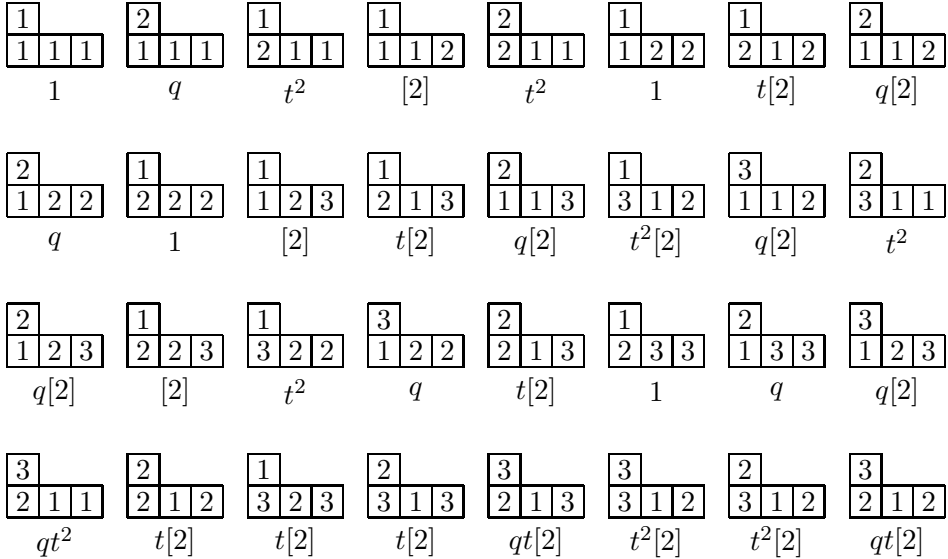
\begin{figure}
\begin{center}
\begin{tikzpicture}[node distance=2cm]
\def \sh {1.6};
\node at (0,0) {$\tableau{ 1\\1&1&1 }$};
\node at (0,-.75) {$1$};
\node at (\sh,0) {$\tableau{ 2\\1&1&1 }$};
\node at (\sh,-.75) {$q$};
\node at (2*\sh,0) {$\tableau{ 1\\2&1&1 }$};
\node at (2*\sh,-.75) {$t^2$};
\node at (3*\sh,0) {$\tableau{ 1\\1&1&2 }$};
\node at (3*\sh,-.75) {$[2]$};
\node at (4*\sh,0) {$\tableau{ 2\\2&1&1 }$};
\node at (4*\sh,-.75) {$t^2$};
\node at (5*\sh,0) {$\tableau{ 1\\1&2&2 }$};
\node at (5*\sh,-.75) {$1$};
\node at (6*\sh,0) {$\tableau{ 1\\2&1&2 }$};
\node at (6*\sh,-.75) {$t[2]$};
\node at (7*\sh,0) {$\tableau{ 2\\1&1&2 }$};
\node at (7*\sh,-.75) {$q[2]$};
\node at (0,-2) {$\tableau{ 2\\1&2&2 }$};
\node at (0,-2.75) {$q$};
\node at (\sh,-2) {$\tableau{ 1\\2&2&2 }$};
\node at (\sh,-2.75) {$1$};
\node at (2*\sh,-2) {$\tableau{ 1\\1&2&3 }$};
\node at (2*\sh,-2.75) {$[2]$};
\node at (3*\sh,-2) {$\tableau{ 1\\2&1&3 }$};
\node at (3*\sh,-2.75) {$t[2]$};
\node at (4*\sh,-2) {$\tableau{ 2\\1&1&3 }$};
\node at (4*\sh,-2.75) {$q[2]$};
\node at (5*\sh,-2) {$\tableau{ 1\\3&1&2 }$};
\node at (5*\sh,-2.75) {$t^2[2]$};
\node at (6*\sh,-2) {$\tableau{ 3\\1&1&2 }$};
\node at (6*\sh,-2.75) {$q[2]$};
\node at (7*\sh,-2) {$\tableau{ 2\\3&1&1 }$};
\node at (7*\sh,-2.75) {$t^2$};

\node at (0,-4) {$\tableau{ 2\\1&2&3 }$};
\node at (0,-4.75) {$q[2]$};
\node at (\sh,-4) {$\tableau{ 1\\2&2&3 }$};
\node at (\sh,-4.75) {$[2]$};
\node at (2*\sh,-4) {$\tableau{ 1\\3&2&2 }$};
\node at (2*\sh,-4.75) {$t^2$};
\node at (3*\sh,-4) {$\tableau{ 3\\1&2&2 }$};
\node at (3*\sh,-4.75) {$q$};
\node at (4*\sh,-4) {$\tableau{ 2\\2&1&3 }$};
\node at (4*\sh,-4.75) {$t[2]$};
\node at (5*\sh,-4) {$\tableau{ 1\\2&3&3 }$};
\node at (5*\sh,-4.75) {$1$};
\node at (6*\sh,-4) {$\tableau{ 2\\1&3&3 }$};
\node at (6*\sh,-4.75) {$q$};
\node at (7*\sh,-4) {$\tableau{ 3\\1&2&3 }$};
\node at (7*\sh,-4.75) {$q[2]$};

\node at (0*\sh,-6) {$\tableau{ 3\\2&1&1 }$};
\node at (0*\sh,-6.75) {$qt^2$};
\node at (1*\sh,-6) {$\tableau{ 2\\2&1&2 }$};
\node at (1*\sh,-6.75) {$t[2]$};
\node at (2*\sh,-6) {$\tableau{ 1\\3&2&3 }$};
\node at (2*\sh,-6.75) {$t[2]$};
\node at (3*\sh,-6) {$\tableau{ 2\\3&1&3 }$};
\node at (3*\sh,-6.75) {$t[2]$};
\node at (4*\sh,-6) {$\tableau{ 3\\2&1&3 }$};
\node at (4*\sh,-6.75) {$qt[2]$};
\node at (5*\sh,-6) {$\tableau{ 3\\3&1&2 }$};
\node at (5*\sh,-6.75) {$t^2[2]$};
\node at (6*\sh,-6) {$\tableau{ 2\\3&1&2 }$};
\node at (6*\sh,-6.75) {$t^2[2]$};
\node at (7*\sh,-6) {$\tableau{ 3\\2&1&2 }$};
\node at (7*\sh,-6.75) {$qt[2]$};

\end{tikzpicture}
\end{center}
        \caption{We compute
		$\widetilde{H}_{2,1,1}(x_1,x_2,x_3;q,t)$ by 
		adding the weights of the sorted tableaux of shape $\lambda=(3,1)$.
		In the figure above, we've listed $t^{\inv(\sigma, \lambda)} q^{\maj, \lambda)} \perm_t(\sigma, \lambda)$ below each tableau, but have omitted the basement and $x^{\sigma}$ to save space.  Here $[i]$ denotes
        $[i]_t$. }\label{fig:alternativeDef}
\end{figure}
\end{example}

To prove Theorem \ref{thm:MacPerm2}, we define \emph{inversion flip operators} which act on fillings of a given shape. These operators fix the $\maj$ statistic and change the $\inv$ statistic by one; in other words, they change the number of counterclockwise inversion triples  by one. Our operators are a generalization of the \emph{inversion flip} move introduced by Loehr and Niese~\cite{LoeNie12} to prove two-column recursions for Macdonald polynomials. Details will appear in the long version of the paper.

\section{A compact formula for integral  Macdonald polynomials}\label{sec:Jcompact}

In this section, we provide a compact formula for the integral form 
Macdonald polynomials $J_{\mu}(X;q;t)$.
We first recall the formula for $J_{\mu}(X;q,t)$ from \cite{HHL04}, with the notation conventions in \cite[Appendix A]{Hag08}. In this formula, for a given filling $\sigma$, the statistic $\coinv(\sigma)$ counts counterclockwise inversion triples of types A and B: 
\begin{itemize}
\item Type A triples consist of the three cells (if they are present in the diagram) $(v,r)$, $(u,r)$, and $(u,r-1)$ for $u<v$. Set $a=\sigma((v,r))$, $b=\sigma((u,r))$, and $c=\sigma((u,r-1))$.
\item Type B triples consist of the three cells (if they are present in the diagram) $(v,r-1)$, $(u,r)$, and $(u,r-1)$ for $u>v$. Set $a=\sigma((v,r-1))$, $b=\sigma((u,r))$, and $c=\sigma((u,r-1))$.
\end{itemize}
We say that a type A or B triple is a \emph{coinversion triple} if one of the following occurs:
\begin{align*}
a < b \leq c\ \ {\rm or}\ \ 
c < a < b \ \ {\rm or}\ \ 
b \leq c < a.
\end{align*}
Note that type A triples coincide with the counterclockwise inversion triples defined in Section \ref{sec:def}.

\begin{theorem}[{\cite[Theorem A.15]{Hag08}}]
	The integral form Macdonald polynomial is given by 
\begin{align} 
\label{Jmu3}
J_{\mu}(X;q,t) =
(1-t)^{\ell (\mu)} \sum_{ \text{nonattacking fillings $\sigma$ 
of $\mu ^{\prime}$ }} x^{\sigma}
q^{\text{maj}(\sigma)} 
t^{\coinv(\sigma)} \\
\notag
\times \prod _{\substack{s \in \mu ^{\prime}\\ \text{$s$ not in row $1$}\\ \sigma (s) = \sigma (\text{South}(s))} }
(1- q^{\text{leg}(s)+1}t^{\text{arm}(s)+1})
\prod _{ \substack{s \in \mu ^{\prime}\\ \text{$s$ not in row $1$}\\ \sigma (s) \ne \sigma ( \text{South}(s) ) } }(1-t),
\end{align}
where $\ell (\mu)$ is the number of parts of $\mu$, and for a cell $s$ not in row $1$,  $\text{South}(s)$ denotes the cell directly below $s$ in the same column as $s$.   By $\text{arm}(s)$ we mean the number of cells in the same row as $s$ and in a column to the right of $s$ whose height is not larger than the column containing $s$, plus (for those $s$ not in the bottom row) the number of cells in the row just below $s$, and in a column to the left of $s$ whose height is strictly smaller than the column containing $s$. Here the sum is over all nonattacking filings of $\mu^{\prime}$ (there is no basement in these fillings). The statistics $\maj$ and $\coinv$ are defined in \cite[Appendix A: pp.~124,~137]{Hag08}
\end{theorem}

In \cite{HHL04} the authors note that the right-hand-side of \eqref{Jmu3} actually yields a correct formula for $J_{\mu}$ if we replace $\mu$ everywhere by $\alpha$, where $\alpha$ is any weak composition of $n$ into $n$ parts satisfying $\text{dec}(\alpha) = \mu$.  (See~\cite[Appendix A]{Hag08} for details about how to extend the notion of $\coinv$ to increasing column heights using Type $B$ triples.)  In fact, the most efficient formula for $J_{\mu}$ seems to be the case where $\alpha = \inc(\mu)$, in which case one can check that \eqref{Jmu3} becomes identical to Lenart's formula for $J_{\mu}$ \cite{Lenart} (which he proved under the additional assumption that $\mu$ has distinct parts).

One unpleasant feature of all these formulas for $J_{\mu}$  is that for the special case $\mu=1^n$, 
\begin{equation}
\label{simple}
J_{1^n}(X;q,t)=x_1x_2\cdots x_n (1-t)(1-t^2)\cdots (1-t^n),
\end{equation}
while the formula \eqref{Jmu3} reduces to
\begin{equation*}
 x_1x_2\cdots x_n (1-t)^n \sum_{\sigma \in S_n}t^{\coinv(\sigma)},
\end{equation*}
a sum of $n!$ terms.  In this section, we show how the identity
\begin{align}
\label{newPmu}
P_{\mu}(X;q,t) =
\sum_{\alpha:\ \text{dec}(\alpha)= \mu} E_{ \inc(\alpha) }^{ \beta (\alpha) }( X;q,t)
\end{align}
yields a corresponding formula for $J_{\mu}$ which, when applied to the case $\mu=1^n$, has only one term---identity \eqref{simple}.
\begin{definition}   
Let $\alpha$ be a weak composition of $n$ into $n$ parts.   We say a nonattacking filling $\sigma$ of $\inc(\alpha)$ (with or without a basement) is \emph{ordered} if in the bottom row of  $\inc(\alpha)$, entries of $\sigma$ below columns of the same height are strictly decreasing when read left to right.   
\end{definition}
Figure \ref{nonatt3b} shows an ordered, nonattacking filling of shape $\alpha=(1,2,2,2,3)$.  The $7$ coinversion triples for the filling are 
$(1,6,7)$, $(3,6,7)$, $(5,6,7)$, $(6,7,9)$, $(1,2,3)$, $(1,2,9)$, and $(3,5,7)$.

\begin{figure}
\begin{center}
\begin{tikzpicture}[scale=.5]
\cell043
\cell147 \cell132 \cell121 \cell 115
\cell246 \cell231 \cell223 \cell215 \cell209
\end{tikzpicture}

\caption{An ordered nonattacking filling with $\maj=3$ and $\coinv= 7$.}
\label{nonatt3b}
\end{center}
\end{figure}

Recall that Macdonald's definition of the integral form $J_{\mu}(X;q,t)$ is
\begin{align}
\label{JDef}
J_{\mu}(X;q,t) &= P_{\mu}(X;q,t) \text{PR}1(\mu),
\end{align}
where 
\begin{equation*}
\label{PR1}
\text{PR}1(\mu) = \prod _{s \in \mu} (1-q^{\text{arm}(s)}t^{\text{leg}(s)+1}) = \prod _{s \in \mu^{\prime}} (1-q^{\text{leg}(s)}t^{\text{arm}(s)+1}).
\end{equation*}

\begin{proposition} 
\label{13}
For $\alpha$ a weak composition,  define
\begin{equation*}
\label{EQ}
\text{PR}2(\alpha) =  \prod_{i\ge 1} (t;t)_{m_i}
\prod_{\substack{s \in \inc(\alpha)^{\prime}\\ \text{$s$ not in the bottom row} }  } (1-q^{\text{leg}(s)+1} t^{\text{arm}(s)+1}),
\end{equation*}
where for $i\ge 1$, $m_i$ is the number of times $i$ occurs in $\alpha$. 
Then if $\mu$ is any partition, $\text{PR}1(\mu) = \text{PR}2(\inc(\mu))$.
\end{proposition}

\begin{definition}
Given a composition $\alpha$ of $n$ into $n$ parts, we define the {\it integral form} version of $E_{\inc(\alpha)} ^{\beta (\alpha)}(X;q,t)$ as
\begin{align}
\label{integralE}
\mathcal {E}_{\inc(\alpha)} ^{\beta (\alpha)}(X;q,t) = PR2(\inc(\alpha))E_{\inc(\alpha)} ^{\beta (\alpha)}(X;q,t).
\end{align}
\end{definition}

Recall the following combinatorial formula 
for $E_{\alpha}^{\sigma}(X;q,t)$ in
\cite{Al16};
\begin{align}
\label{formula}
E_{\inc(\alpha)}^{\beta (\alpha)}(X;q,t) = \sum_{\sigma} x^{\sigma} \wt(\sigma),
\end{align}
where the sum is over all nonattacking fillings $\sigma$ of the diagram of $\alpha$ 
whose $i$th column is of height $\inc(\alpha )_i$, with basement $\beta (\alpha)$.  The weight $\wt(\sigma)$ is
\begin{align}
\label{weight}
\wt(\sigma) = q^{ \text{maj} (\sigma)} t^{ \coinv(\sigma) }\prod_{s:\ \sigma(s) \ne \sigma(\text{South}(s))} \frac {1-t}{1-q^{\text{leg(s)}+1}t^{\text{arm(s)}+1}}.
\end{align}

It follows from the above formula that $\mathcal {E}_{\inc(\alpha)} ^{\beta (\alpha)}(X;q,t)$ 
is $\prod_{i} (t;t)_{m_i}$ times an element of $\mathbb Z [x_1,\ldots ,x_n,q,t]$.  To see this, note that every nonattacking filling $\sigma$ of $\inc(\alpha)$ has the property that each entry in the bottom row is equal to the entry in the basement directly below it, and hence doesn't contribute anything to the product in \eqref{weight} defining $\wt(\sigma)$, while if any entry above the bottom row satisfies $\sigma (s) \not= \sigma (\text{South}(s))$, then the associated  factor  
$(1-q^{\text{leg}(s)+1}t^{\text{arm}(s)+1})$ is the exact term in the coefficient of
 $E_{\inc(\alpha)} ^{\beta (\alpha)}(X;q,t)$ from \eqref{integralE} above corresponding to cell $s$. In fact this argument shows that
\begin{align}
\label{IFE}
\mathcal {E}_{\inc(\alpha)} ^{\beta (\alpha)}(X;q,t) = 
\prod_{i}(t;t)_{m_i}  
\sum_{\substack{\text{ordered, nonattacking fillings $\sigma$ 
of inc$(\alpha) ^{\prime}$ }\\ \text{basement $\beta (\alpha)$} } } 
x^{\sigma}q^{\text{maj}(\sigma) } t^{\coinv(\sigma )} \\
\notag
\times \prod _{\substack{s \in \inc(\alpha) ^{\prime},\ \text{$s$ not in row $1$} \\ 
\sigma (s) = \sigma (\text{South}(s))} }
(1- q^{\text{leg}(s)+1}t^{\text{arm}(s)+1})
\prod _{ \substack{s \in \inc(\alpha)^{\prime},\ \text{$s$ not in row $1$}\\ \sigma (s) \ne \sigma ( \text{South}(s) ) } }(1-t),
\end{align}

\begin{corollary} The formula for $J_{\mu}$ has the following more compact version:
\begin{align} 
\label{Jmu4}
J_{\mu}(X;q,t) =
\prod_{i}(t;t)_{m_i}  
\sum_{ \text{ ordered, nonattacking fillings $\sigma$ 
of inc$(\mu) ^{\prime}$ } } 
x^{\sigma}q^{\text{maj}(\sigma)} t^{\coinv(\sigma)} \\
\notag
\times \prod _{\substack{s \in \mu ^{\prime},\ \text{$s$ not in row $1$}\\ \sigma (s) = \sigma (\text{South}(s))} }
(1- q^{\text{leg}(s)+1}t^{\text{arm}(s)+1})
\prod _{ \substack{s \in \mu ^{\prime},\ \text{$s$ not in row $1$}\\ \sigma (s) \ne \sigma ( \text{South}(s) ) } }(1-t).
\end{align}
\end{corollary}
Note that \eqref{Jmu4} implies the (as far as we know) new fact that $J_{\mu}(X;q,t)$ is $\prod_{i}(t;t)_{m_i}$ times an element of
$\mathbb Z [x_1,\ldots, x_n,q,t]$.
\begin{proof}
Start by multiplying both sides of \eqref{newPmu} by $\text{PR}2(\inc(\mu))$.  The left-hand-side then becomes $J_{\mu}(X;q,t)$ by \eqref{JDef} and Proposition \ref{13}.   
The summand on the right hand side becomes $\mathcal {E}_{\inc(\alpha)} ^{\beta (\alpha)}(X;q,t)$, which by \eqref{IFE} equals the portion of \eqref{Jmu4} which has bottom row determined by $\beta(\alpha)$.
\end{proof}

\section{A quasisymmetric Macdonald polynomial}\label{sec:qsym}

Recall that the ring of \emph{quasisymmetric functions} is a graded ring 
which  contains within it 
the ring of symmetric functions.  
The ring of quasisymmetric functions
 has multiple
distinguished bases, indexed by (strong) compositions, 
one of which is the \emph{quasisymmetric Schur
functions} 
$\text{QS}_{\gamma}(X)$ 
introduced by the second and fourth authors, together with Luoto and van Willigenburg 
\cite{HLMV09}.  
The authors showed that $\text{QS}_{\gamma}(X)$ is quasisymmetric,
and that each (symmetric) Schur function $s_{\lambda}(X)$ 
is a positive sum of quasisymmetric 
Schur functions.  In light of this, and the fact that Macdonald 
polynomials expand positively in terms of Schur polynomials, it is 
natural to ask if there is a notion of a 
\emph{Macdonald quasisymmetric polynomial} $G_{\gamma}(X; q, t)$ such that:
\begin{enumerate}
	\item[(A)]\label{a} The symmetric Macdonald polynomial $P_{\lambda}(X; q, t)$ is a positive
	     sum of Macdonald quasisymmetric polynomials; 
     \item[(B)]\label{b} $G_{\gamma}(X; q, t)$ is quasisymmetric; 
     \item[(C)]\label{c} $G_{\gamma}(X; 0, 0)$ is the quasisymmetric Schur function 
		$\text{QS}_{\gamma}(X)$;
 \item[(D)] \label{d} $G_{\gamma}(X; q, t)$ has a combinatorial formula
along the lines of the ``HHL'' formula for the $E_{\alpha}$ \cite{HHL08}, or its compact version from \cite{CMW18}.
\end{enumerate}
We show in this section that the answer to this question is yes.

Given a permutation $\sigma \in S_n$,  let $E_{\alpha}^{\sigma}(X;q,t)$ be the
 permuted-basement nonsymmetric Macdonald polynomial defined in \cite{thesis} and studied in  \cite{Al16,CMW18}, and let
 $F_{\alpha}(X;q,t) = E_{\inc(\alpha)}^{\beta (\alpha)}(X;q,t)$.  For any partition $\lambda$ of $n$, from \cite{CMW18} we have that 
 \begin{align}
 \label{basic}
 P_{\lambda}(X;q,t)  = \sum_{\alpha :\ \text{dec}(\alpha)= \lambda} F_{\alpha}(X;q,t),
 \end{align}
where 
the sum is over all weak compositions $\alpha$ whose positive parts are a 
rearrangement of the parts of $\lambda$. 

Note that if $\text{id}=(1,2,\ldots ,n)$ and  $w_0=(n,n-1,\ldots,1)$ are the identity permutation and permutation of maximal length in $S_n$,
respectively, then
$E_{\alpha}^{ \text{id} }(X;0,0)$ is the Demazure atom and $E_{\alpha}^{ w_0 }(X;0,0)$ the Demazure character. (In the common notation for Demazure characters, 
i.e., key polynomials, one reverses the vector $\alpha$, i.e.,  the key polynomial corresponding to $\alpha$ would be $E_{(\alpha _n, \ldots ,\alpha _1)}^{w_0}(X;0,0)$.)

Motivated by \eqref{basic}, we have the following definition and theorem.

\begin{theorem}
	We define the \emph{quasisymmetric Macdonald polynomial}
	$G_{\gamma}(X;q,t)$  to be
\begin{align}
\label{Gdef}
G_{\gamma}(X;q,t) &= \sum_{\alpha:\ \alpha ^{+}= \gamma}  F_{\alpha}(X;q,t) \\
&=  \sum_{\alpha:\ \alpha ^{+}= \gamma}  E_{\inc(\alpha)}^{\beta (\alpha )}(X;q,t),\nonumber
\end{align}
where the sum is over all weak compositions $\alpha$ for which $\alpha^{+}=\gamma$.
	Then $G_{\gamma}(X; q, t)$ satisfies properties 
	\rm{(A), (B), (C), and (D)}.
\end{theorem}

The fact that $G_{\gamma}(X; q, t)$ satisfies \rm{(A)} follows from \eqref{basic}. There are several combinatorial proofs that $G_{\gamma}(X; q,t)$ is quasisymmetric and hence satisfies \rm{(B)}. One proof uses \eqref{formula} and a notion of \emph{packed} nonattacking fillings. Another proof uses the multiline queues from \cite{CMW18}. To see that  $G_{\gamma}(X; q, t)$ satisfies \rm{(C)}, recall that for $\gamma$ a strong composition of $n$, $\text{QS}_{\gamma}$ is defined by the equation
\begin{align}
\label{QS}
\text{QS}_{\gamma}(X) = \sum_{\alpha:\ \alpha ^{+}= \gamma} E_{\alpha}^{\text{id}}(X;0,0).
\end{align}
So to verify \rm{(C)}, it suffices to show that $F_{\alpha}(X; 0, 0) = E_{\alpha}^{\text{id}}(X; 0, 0)$.  We actually show the stronger statement that 
\begin{equation} \label{tatom}
	F_{\alpha}(X; 0, t) = E_{\alpha}^{\text{id}}(X; 0, t),
\end{equation}
	where $E_{\alpha}^{\text{id}}(X; 0, t)$ is the Demazure $t$-atom.
To prove this, one can use induction, together with the action of the Hecke operators.
In particular, by \cite{CMW}, $T_i F_{\alpha} = F_{s_i \alpha}$ if $\alpha_i >\alpha_{i+1}$; compare with 
\cite[Corollary 26]{Al16}.

It would be interesting to find a connection between the quasisymmetric Macdonald polynomials
$G_{\gamma}(X; q, t)$ that we introduce in this paper, and other objects in the literature.
We note that 
our $G_{\gamma}$ are different from the duals of the noncommutative symmetric function analogues of Macdonald polynomials introduced in~\cite{BerZab05}; 
we also do not see a connection to the noncommutative Hall-Littlewood polynomials studied in~\cite{Hiv00}.


\end{document}